\documentclass[12pt]{amsart}
\usepackage{amssymb, latexsym, amsthm, amsmath, amsxtra}

\theoremstyle{plain}

\newtheorem{thm}{Theorem}        
\newtheorem{theorem}[thm]{Theorem}
\newtheorem{lem}[thm]{Lemma}

\theoremstyle{definition}
\newtheorem{prob}{Problem}

\title[Laplace-Beltrami operator in adiabatic
limit]{The spectrum of the Laplace-Beltrami operator on the
two-dimensional torus in adiabatic limit}
\author{Andrey A. Yakovlev}
\address{Department of Mathematics,\\ Ufa State Aviation Technical
University, 12 K.Marx str., Ufa, 450000, Russia}
\email{yakovlevandrey@yandex.ru}
\thanks{Supported by Russian Foundation of Basic Research
(grant no. 06-01-00208)}
\date{}

\begin{document}

\begin{abstract}
We obtain an asymptotic formula for the eigenvalue distribution
function of the Laplace-Beltrami operator on the two-dimensional
torus in the adiabatic limit given by a Kronecker foliation.
Related problems in number theory are discussed.
\end{abstract}

\maketitle

\section*{Introduction}
Consider the two-dimensional torus
$\mathbb{T}^2=\mathbb{R}^2/\mathbb{Z}^2$ with the coordinates
$(x,y)\in \mathbb{R}^2$, considered modulo integer translations,
and the Euclidean metric $g$ on $\mathbb{T}^2$:
\[
g=d x^2+d y^2.
\]

Let $\widetilde{X}$ be the vector field on $\mathbb{R}^2$ given by
\[
\widetilde{X}=\frac{\partial}{\partial x}+\alpha
\frac{\partial}{\partial y},
\]
where $\alpha\in \mathbb R$. Since $\widetilde{X}$ is translation
invariant, it determines a vector field $X$ on $\mathbb{T}^2$.
The vector field $X$ defines a foliation $\mathcal F$ on
$\mathbb{T}^2$. The leaves of $\mathcal F$ are the images of the
parallel lines $\widetilde{L}=\{(x_{0}+t, y_{0}+t\alpha)
:t\in\mathbb{R},(x_0,y_0)\in \mathbb R^2\}$ with the slope
$\alpha$ under the projection $\mathbb{R}^2 \rightarrow
\mathbb{T}^2$.

In the case when $\alpha$ is rational, all leaves of $\mathcal F$
are closed and are circles, and $\mathcal F$ is given by the
fibers of a fibration of $\mathbb{T}^2$ over $\mathbb{S}^1$. In
the case when $\alpha$ is irrational, all leaves of $\mathcal F$
are everywhere dense in $\mathbb{T}^2$.

Let $F=T{\mathcal F}$ be the tangent bundle to $\mathcal F$:
\[
F_{(x,y)}=\{ (t,\alpha t):t\in \mathbb{R} \}=<(1,\alpha)>\subset
 T_{(x,y)}\mathbb{T}^2,
\]
and $H=F^{\bot}=<(-\alpha,1)>$ is the orthogonal complement for
$F$. Thus, we have the decomposition of the bundle
$T\mathbb{T}^2$ into direct sum $T\mathbb{T}^2=F\bigoplus H$
\[
T_{x,y}\mathbb{T}^2\ni(X,Y)=\frac{X+\alpha
Y}{1+\alpha^2}(1,\alpha)+\frac{-\alpha X+Y}{1+\alpha^2}(-\alpha,1)
\]
and the corresponding decomposition of the Riemannian metric
$g=g_{F}+g_{H}$. Define a one-parameter family $g_{h}$ of
Riemannian metrics on $\mathbb{T}^2$ by the formula
\begin{equation*}
g_{h}=g_{F} + {h}^{-2}g_{H}, \quad 0 < h \leq 1.
\end{equation*}
In the coordinates $(x,y)$, the Riemannian metric $g_h$ is given
by
\[
g_h=
  \frac{1+h^{-2}\alpha^2}{1+\alpha^2} dx^2 + 2 \alpha
  \frac{1-h^{-2}}{1+\alpha^2}dx dy+
  \frac{\alpha+h^{-2}}{1+\alpha^2}dy^2.
\]

Recall that, for an arbitrary $n$-dimensional manifold $M$ with a
Riemannian metric $g$, an elliptic, second order, differential
operator $\Delta=\Delta_g$ in the space $C^\infty_0(M)$, which is
called the Laplace-Beltrami operator, is defined. In local
coordinates, this operator is given by the formula
\[
\Delta=-\frac{1}{\sqrt{g}}\sum_{i,j=1}^n\frac{\partial}{\partial
x^i}\left(g^{ij}\sqrt{g}\frac{\partial}{\partial x^j}\right),
\]
where $(g_{ij})$ is the matrix of the Riemannian metric, $g$ is
the determinant of the matrix $(g_{ij})$, $(g^{ij})$ is the
inverse of the matrix $(g_{ij})$.

In our case, the Laplace operator, defined by the metric $g_h$,
has the form
 \[
\Delta_h=-\frac{1}{1+\alpha^2}\left( \frac{\partial}{\partial
x}+\alpha\frac{\partial}{\partial y}
\right)^2-\frac{h^{2}}{1+\alpha^2}\left(
-\alpha\frac{\partial}{\partial x}+\frac{\partial}{\partial y}
\right)^2=\Delta_F+h^2\Delta_H,
 \]
where $\Delta_F$ is the tangential Laplace operator:
 \[
\Delta_F=-\frac{1}{1+\alpha^2}\left( \frac{\partial}{\partial
x}+\alpha\frac{\partial}{\partial y} \right)^2,
\]
and $\Delta_H$ is the transverse Laplace operator:
\[
\Delta_H=-\frac{h^{2}}{1+\alpha^2}\left(
-\alpha\frac{\partial}{\partial x}+\frac{\partial}{\partial y}
\right)^2.
 \]

The main result of the paper is the calculation of the
asymptotics of the spectrum distribution function $N_h(\lambda)$
of the operator $\Delta_h$ for a fixed $\lambda\in\mathbb{R}$ and
for $h\rightarrow 0$.

\begin{theorem}\label{th:main}
The following asymptotic formula for the spectrum distribution
function $N_h(\lambda)$ of the operator $\Delta_h$ for a fixed
$\lambda\in\mathbb{R}$ holds:
\medskip\par
1. For $\alpha\not\in\mathbb{Q},$
\[
N_{h}(\lambda ) =\frac{1}{4\pi}h^{-1}\lambda+o(h^{-1}), \quad
h\rightarrow 0.
\]

2. For $\alpha\in\mathbb{Q}$ of the form $\alpha=\frac{p}{q}$,
where $\operatorname{G.C.D.}(p,q)=1$,
\[
 N_{h}(\lambda ) =h^{-1}
\sum_{k:\frac{4\pi^2}{p^2+q^2}k^2<\lambda}
\frac{1}{\pi\sqrt{p^2+q^2} }(\lambda -
\frac{4\pi^2}{p^2+q^2}k^2)^{1/2}+o(h^{-1}), \quad h\rightarrow 0.
\]
\end{theorem}

We will give two proofs of Theorem~\ref{th:main}. The first proof
is based on a general result on the asymptotic behavior of the
spectrum distribution function of the Laplace operator on a
compact manifold equipped with a Riemannian foliation in the
adiabatic limit, obtained in the work \cite{adiab}. The second
proof is considerably simpler and uses only elementary facts of
analysis.

Regarding various notions and facts, connected with the foliation
theory see, for instance, the survey paper \cite{survey} and the
bibliography cited there.

The results of the paper are announced in \cite{1}.

The author is grateful to Kordyukov Yu.A. for the setting of the
problem and useful discussions.

\section{Adiabatic limits for Riemannian foliations}
Let $(M,{\mathcal F})$ be a closed foliated manifold, $\dim M =
n$, $\dim {\mathcal F} = p$, $p+q=n$, equipped with a Riemannian
metric $g_M$. Let $F=T{\mathcal F}$ be the integrable distribution
of tangent $p$-planes to ${\mathcal F}$ in $TM$, and $H=F^{\bot}$
the orthogonal complement to $F$. Thus, we have a direct sum
decomposition of $TM$:
\begin{equation*}
TM=F\oplus H.
\end{equation*}
This decomposition induces a decomposition of the metric $
g_{M}=g_{F}+g_{H}$. Define a one-parameter family
$g_{h}$ of metrics on $M$ by the formula
\begin{equation*}
g_{h}=g_{F} + {h}^{-2}g_{H}, \quad 0 < h \leq 1.
\end{equation*}

For any $h>0$ we consider the Laplace-Beltrami operators on
functions, determined by the metric $g_h$:
\begin{equation*}
\Delta_{h}=d^{*}_{g_h}d,
\end{equation*}
where $d :C^{\infty}(M)\rightarrow C^{\infty}(M,\Lambda^1 T^{*}M)$
is the de Rham differential, $d^{*}_{g_h}$ is the adjoint operator
of $d$ with respect to the inner products, defined by the metric
$g_{h}$. The operator $\Delta_{h}$ is a self-adjoint, elliptic,
differential operator with the positive definite, scalar principal
symbol in the Hilbert space $L^2(M,g_h)$ of square integrable
differential forms on $M$, endowed with the inner product, induced
by the metric $g_h$. According to the standard perturbation theory
there exists a countable family of analytic functions
$\lambda_i(h)$ such that, for any $h>0$,
\begin{equation*}
{\rm spec}\,\Delta_h =\{\lambda_i(h):i=0,1,\ldots\},
\end{equation*}
taking into account multiplicities. Introduce the eigenvalue
distribution function of the operator $\Delta_h$:
\begin{equation*}
N_h(\lambda)=\sharp \{i:\lambda_i(h)\leq \lambda\}
\end{equation*}

Suppose that the foliation $\mathcal F$ is Riemannian and $g_M$ is
a bundle-like Riemannian metric. Let $dx$ be the Riemannian volume
form on $M$.

Denote by $\Delta_F$ the tangential Laplace operator on
$C^{\infty}(M)$ and by $\Delta_L$ its restriction to the leaf $L$
of $\mathcal F$. The operator $\Delta_L$ coincides with the
Laplace operator on $L$, defined by the induced metric of the
leaf.

Let $G$ be the holonomy groupoid of $\mathcal F$. Let $\lambda_L$
denote the Riemannian volume form on a leaf $L$ given by the
induced metric, and $\lambda^x$ denote its lift to the holonomy
covering $G^x$ for any $x\in M$.

Denote by $e_\lambda(\gamma), \gamma\in G, \lambda \in
{\mathbb{R}}$ the tangential kernel of the spectral projection of
the tangentially elliptic operator $\Delta_{F}$, corresponding to
the semi-axis $(-\infty ,\lambda )$. For every $x\in M$ the
function $e_\lambda(\gamma^{-1}\gamma'), (\gamma, \gamma')\in G^x$
coincides with the kernel of the spectral projection of the
operator $\Delta_{G^x}$, which is the lift of $\Delta_L$ to the
holonomy covering $G^x$, relative to the Riemannian volume form
$\lambda^x$. Let $e_\lambda(x), x\in  M, \lambda \in {\mathbb{R}}$
be the restriction of $e_\lambda(\gamma)$ to $M$. For any $\lambda
\in {\mathbb{R}}$, the function $e_\lambda$ is a bounded
measurable function on $M$.

The spectrum distribution function $N_{\mathcal F}(\lambda )$ of
$\Delta_{F}$ is defined by the formula
\[
N_{\mathcal F}(\lambda )=\int_M e_\lambda(x)\,dx.
\]

\begin{theorem}
\label{ad:intr} Let $(M,{\mathcal F})$ be a Riemannian foliation,
equipped with a bundle-like metric $g_M$. The following asymptotic
formula holds:
\begin{equation}
\label{ad:eig1} N_{h}(\lambda ) =h^{-q}
\frac{(4\pi)^{-q/2}}{\Gamma((q/2)+1)} \int_{-\infty}^{\lambda}
(\lambda - \tau )^{q/2}d_{\tau}N_{\mathcal F}(\tau )+o(h^{-q}),
\quad h\rightarrow 0.
\end{equation}
\end{theorem}

\section{Proof of Theorem~\ref{th:main}}
\subsection{The first proof} In this section we derive the asymptotic formula for
$N_h(\lambda)$ on the torus from Theorem~\ref{ad:intr}.

The induced metric on the leaf $x=x_0+t,y=y_0+\alpha t$ through
a point $(x_0,y_0)$ has the form
$g_\mathcal{F}=(1+\alpha^2)dt^2$. Thus, the leafwise Laplace operator is given by
\[
\Delta_F=-\frac{1}{1+\alpha^2}\left(\frac{\partial}{\partial
x}+\alpha \frac{\partial}{\partial y}\right)^2.
\]

1. $\alpha\not\in\mathbb{Q}.$ In this case $G=\mathbb{T}^2\times
\mathbb{R}$. The source and the range maps $s,r:G\rightarrow
\mathbb{T}^2$ are defined for any $\gamma=(x,y,t)\in G$ by
$s(\gamma)=(x-t,y-\alpha t)$ and $r(\gamma)=(x,y)$. For any
$(x,y)\in \mathbb{T}^2$, the set $G^{(x,y)}$ coincides with the
leaf $L_{(x,y)}$ through $(x,y)$ and is diffeomorphic to $\mathbb
R$:
\[
L_{(x,y)}=\{(x-t,y-\alpha t): t\in {\mathbb R}\}.
\]
The Riemannian volume form $\lambda^{(x,y)}$ on $L_{(x,y)}$ equals
$\sqrt{1+\alpha^2}\,dt$. Finally, the restriction of the operator
$\Delta_F$ to each leaf $L_{(x,y)}$ coincides with the operator
\[
A=-\frac{1}{1+\alpha^2}\frac{d^2}{dt^2},
\]
acting in the space $L^2({\mathbb{R}}, \sqrt{1+\alpha^2}\,dt)$.

Let us find the spectral projections of $A$ in $L^2(\mathbb{R},
\sqrt{1+\alpha^2}\,dt)$. The following equality is valid:
\[
A_1U(\xi)=FAF^{-1}U(\xi)=\frac{|\xi|^2}{1+\alpha^2}U(\xi), \quad
U\in C^\infty_0(\mathbb{R}), \] where $F$ is the Fourier
transform. A similar equality holds for the spectral projections:
\[
\chi_\lambda(A_1)=F\chi_\lambda(A)F^{-1},
\]
where $\chi_\lambda$ is the characteristic function of the semi-axis
$(-\infty,\lambda)$. It is easy to see that the operator
$\chi_\lambda(A_1)$ is the multiplication operator by the the function
\[
\chi_\lambda\left(\frac{|\xi|^2}{1+\alpha^2}\right)=\begin{cases}
                            1,\ \text{if}\ \frac{|\xi|^2}{1+\alpha^2}<\lambda,\\
                            0,\ \text{if}\ \frac{|\xi|^2}{1+\alpha^2}\geq\lambda.
                         \end{cases}
\]
Thus, we get that the spectral projection $\chi_\lambda(A)$ of the
operator $A$ is of the form:
\begin{align*}
\chi_\lambda(A)U(t)&=F^{-1}\chi_\lambda\left(\frac{|\xi|^2}{1+\alpha^2}\right)FU(t)\\
&=
\frac{1}{2\pi}\int_{\mathbb{R}^2}{\exp{[i(t-t_1)\xi]}\chi_\lambda\left(\frac{|\xi|^2}{1+\alpha^2}\right)U(t)dt_1
d\xi}.
\end{align*}
The kernel of the spectral projection $\chi_\lambda(A)$ of $A$ in the space $L^2(\mathbb{R})$ (relative to the Riemannian volume form  $\sqrt{1+\alpha^2}\,dt$) is given by the formula
\[
E_\lambda(t,t_1)=\frac{1}{2\pi\sqrt{1+\alpha^2}}\int_{\mathbb{R}}{\exp{[i(t-t_1)\xi]}
\chi_\lambda\left(\frac{|\xi|^2}{1+\alpha^2}\right)d\xi}.
\]
The tangential kernel $e_\lambda$ of the operator $\Delta_F$ is
related with $E_\lambda$ as follows. For any $\gamma=(x,y,t)\in
G=\mathbb{T}^2\times \mathbb{R}$,
\[
e_\lambda(\gamma)=E_{\lambda}(0,t)=\frac{1}{2\pi\sqrt{1+\alpha^2}}
\int_{\mathbb{R}}{\exp{(-it\xi)}\chi_\lambda\left(\frac{|\xi|^2}{1+\alpha^2}\right)d\xi}.
\]
The restriction of the tangential kernel $e_\lambda$ to $\mathbb{T}^2$ is given by
\[
e_\lambda(x,y)=E_{\lambda}(0,0)=\frac{1}{2\pi\sqrt{1+\alpha^2}}\int_{\mathbb{R}}\chi_\lambda\left(\frac{|\xi|^2}{1+\alpha^2}\right)d\xi
= \frac{1}{\pi}\sqrt{\lambda}, \quad \lambda>0.
\]
We get that the spectrum distribution function $N_{\mathcal
F}(\lambda )$ of operator $\Delta_{F}$ has the form:
\[
N_\mathcal{F}(\lambda)=\int_{\mathbb{T}^2}{e_\lambda(x,y)dxdy}=
 \frac{1}{\pi}\sqrt{\lambda}, \quad \lambda>0.
 \]
By Theorem~\ref{ad:intr}, we obtain
\[
N_{h}(\lambda ) =h^{-1} \frac{1}{\pi} \int_{-\infty}^{\lambda}
(\lambda - \tau )^{1/2}d_{\tau}N_{\mathcal F}(\tau )+o(h^{-1}),
\quad h\rightarrow 0.
\]
Hence,
\[
N_{h}(\lambda ) =\frac{1}{4\pi}h^{-1}\lambda+o(h^{-1}), \quad
h\rightarrow 0.
\]

2. $\alpha\in\mathbb{Q}$. Let us consider $\alpha=\frac{p}{q}$,
where $p$ and $q$ are coprime. In this case, the holonomy groupoid
is $\mathbb{T}^2\times (\mathbb{R}/{q\mathbb{Z}})$. The leaf
through an arbitrary point $(x_0,y_0)$ is the circle
$\{(x_0+t,y_0+\alpha t): t\in \mathbb{R}/{q\mathbb{Z}}\}$ of
length $l=\sqrt{p^2+q^2}$.

The restriction of the operator $\Delta_F$ to each $L_{(x,y)}$
coincides with the operator
\[
A=-\frac{1}{1+\alpha^2}\frac{d^2}{dt^2},
\]
acting in the space $L^2(\mathbb{R}/q\mathbb{Z},
\sqrt{1+\alpha^2}\,dt)$.

This operator has a complete orthogonal system of eigenfunctions
$U_j(t)$ the corresponding eigenvalues $\lambda_j $:
\[
\lambda_j=\frac{1}{1+\alpha^2}\frac{4\pi^2}{q^2}j^2=\frac{4\pi^2}{p^2+q^2}j^2,
\quad U_j(t)=\exp [\frac{2\pi i}{q}jt], \quad j\in \mathbb{Z}.
\]
Any function $U(t)\in L^2(\mathbb{R}/q\mathbb{Z},
\sqrt{1+\alpha^2}\,dt)$ can be represented as the sum of a series:
\[
    U(t)=\frac{1}{q}\sum_{k=-\infty}^\infty{\left(\int_0^q{U(t_1)
    \exp[-\frac{2\pi i}{q}k t_1]dt_1}
    \right)\exp[\frac{2\pi i}{q}k t]}
\]
Therefore, the spectral projection $\chi_\lambda(A)$ of
$A$ is given by the formula:
\[
    \chi_\lambda(A) U(t)=\frac{1}{q}\sum_{k:\frac{4\pi^2}{p^2+q^2}k^2<\lambda}
    {\left(\int_0^q{U(t_1) \exp[-\frac{2\pi i}{q}k t_1]dt_1}
    \right)\exp[\frac{2\pi i}{q}k t]}
\]
The kernel of the projection $\chi_\lambda(A)$ in
$L^2(\mathbb{R}/q\mathbb{Z},\sqrt{1+\alpha^2}\,dt)$ is given by
the formula
\[
    E_{\lambda}(t,t_1)=\frac{1}{\sqrt{p^2+q^2}}\sum_{k:\frac{4\pi^2}{p^2+q^2}k^2<\lambda}{\exp[\frac{2\pi
i}{q}k(t-t_1)]}.
    \]
The tangential kernel $e_\lambda$ of $\Delta_F$ is related with
$E_\lambda$ in the following way. For any $\gamma=(x,y,t)\in
G=\mathbb{T}^2\times(\mathbb{R}/q\mathbb{Z})$,
\[
e_\lambda(\gamma)=E_\lambda(0,t)=\frac{1}{\sqrt{p^2+q^2}}\sum_{k:\frac{4\pi^2}{p^2+q^2}k^2<\lambda}{\exp[-\frac{2\pi
i}{q}k t)]}.
\]
The restriction $e_\lambda$ to $\mathbb T^2$ has the form:
\[    e_\lambda(x,y)=E_{\lambda}(0,0)=\frac{1}{\sqrt{p^2+q^2}}\#\{k:\frac{4\pi^2}{p^2+q^2}k^2<\lambda\}.
\]

We get that, in the case of a rational $\alpha$, the spectrum distribution function  $N_{\mathcal F}(\lambda )$ of
$\Delta_{F}$ is of the form:
\[
N_{\mathcal{F}}(\lambda )= \int_{\mathbb{T}^2}{e_\lambda(x,y)dxdy}
=\frac{1}{\sqrt{p^2+q^2}}\#\{k:\frac{4\pi^2}{p^2+q^2}k^2<\lambda\}.
\]
By Theorem~\ref{ad:intr}, we obtain for $h\rightarrow 0$
\[
\begin{split}
 N_{h}(\lambda )
& =h^{-1} \frac{1}{\pi} \int_{-\infty}^{\lambda} (\lambda - \tau
)^{1/2}d_{\tau}N_{\mathcal F}(\tau )+o(h^{-1})\\ & =h^{-1}
\frac{1}{\pi\sqrt{p^2+q^2} }
\sum_{k:\frac{4\pi^2}{p^2+q^2}k^2<\lambda} (\lambda -
\frac{4\pi^2}{p^2+q^2}k^2)^{1/2}+o(h^{-1}).
\end{split}
\]

\subsection{The second proof}

The operator $\Delta_h$ has a complete orthogonal system of eigenfunctions
\[
u_{kl}(x,y)=e^{2\pi i(kx+ly)}, \quad (x,y)\in {\mathbb{T}^2},
\]
with the corresponding eigenvalues
\[
\lambda_{kl}= (2\pi)^2 \left(\frac{1}{1+\alpha^2}(k+\alpha
l)^2+\frac{h^2}{1+\alpha^2}(-\alpha k+l)^2\right), \quad k,l\in
{\mathbb Z}^2.
\]

The eigenvalue distribution function of $\Delta_h$ has the form
\begin{multline*}
N_h= \# \{(k,l)\in {\mathbb Z}^2 :\\ (2\pi)^2
\left(\frac{1}{1+\alpha^2}(k+\alpha
l)^2+\frac{h^2}{1+\alpha^2}(-\alpha k+l)^2\right) < \lambda\}.
 \end{multline*}
Thus we come to the following problem of number theory:
\begin{prob}
Find the asymptotics for $h\rightarrow 0$ of the number of integer
points in the ellipse
\[
\{(\xi,\eta)\in \mathbb{R}^2: (2\pi)^2
\left(\frac{1}{1+\alpha^2}(\xi+\alpha
\eta)^2+\frac{h^2}{1+\alpha^2}(-\alpha
\xi+\eta)^2\right)<\lambda\}.
\]
\end{prob}

Let us introduce new coordinates
\[
\left\{
\begin{matrix}
  \xi^\prime=\frac{1}{\sqrt{1+\alpha^2}}\xi+\frac{\alpha}{\sqrt{1+\alpha^2}}\eta ,\\
  \eta^\prime=-\frac{\alpha}{\sqrt{1+\alpha^2}}\xi+\frac{1}{\sqrt{1+\alpha^2}}\eta . \\
\end{matrix}
\right.
\]
In the new coordinates we will consider a point
$(\xi^\prime,\eta^\prime)$ to be integer, if it has the form
\[
\left\{
\begin{matrix}
  \xi^\prime=\frac{1}{\sqrt{1+\alpha^2}}n+\frac{\alpha}{\sqrt{1+\alpha^2}}m ,\\
  \eta^\prime=-\frac{\alpha}{\sqrt{1+\alpha^2}}n+\frac{1}{\sqrt{1+\alpha^2}}m ,\\
\end{matrix}
\right.
\]
where $(n,m)\in \mathbb{Z}^2$. Hence, the number of such integer
points is determined by the number of all pairs $(n,m)\in
\mathbb{Z}^2$, satisfying the following conditions:
\[
\left\{
\begin{matrix}
  |\frac{1}{\sqrt{1+\alpha^2}}n+\frac{\alpha}{\sqrt{1+\alpha^2}}m|<\frac{\sqrt{\lambda}}{2\pi} ,\\
  |-\frac{\alpha}{\sqrt{1+\alpha^2}}n+\frac{1}{\sqrt{1+\alpha^2}}m|
  < h^{-1}\sqrt{\frac{\lambda}{(2\pi)^2}-(\frac{1}{\sqrt{1+\alpha^2}}n+\frac{\alpha}{\sqrt{1+\alpha^2}}m)^2}. \\
\end{matrix}
\right.
\]
It is clear that, as $h\rightarrow 0$, for
$\alpha\in\mathbb{R}/\mathbb{Q}$ the projections of integer points to the
$\xi^\prime$-axes form an everywhere dense set, but for
$\alpha\in\mathbb{Q}$ a discrete one.

Suppose that $\alpha\in\mathbb{Q}$. Without loss of generality, we will assume  $\alpha=\frac{p}{q}$, where $\text{G.C.D.}(p,q)=1$. By straightforward computations, it is easy to show that a point is integer, if has the form:
\[
(\frac{k}{\sqrt{p^2+q^2}},\frac{-pn_0+qm_0-(p^2+q^2)l}{\sqrt{p^2+q^2}}),
\]
where $(n_0,m_0)$ is an integer point determined by $k$,
$l\in\mathbb{Z}$ is responsible for the number of integer point over
$\xi^\prime$. Since we are looking for the number of integer points in the ellipse, we get the following conditions on $l$:
\begin{multline*}
-h^{-1}\sqrt{\frac{\lambda}{4\pi^2}-\frac{k^2}{p^2+q^2}}+
\frac{-pn_0+qm_0}{\sqrt{p^2+q^2}}<\sqrt{p^2+q^2}l \\ <
h^{-1}\sqrt{\frac{\lambda}{4\pi^2}-\frac{k^2}{p^2+q^2}}+
\frac{-pn_0+qm_0}{\sqrt{p^2+q^2}}.
\end{multline*}
Thus, the number of integer points over $\xi^\prime$ is equal to
\[
2\frac{h^{-1}}{\sqrt{p^2+q^2}}\sqrt{\frac{\lambda}{4\pi^2}-\frac{k^2}{p^2+q^2}}.
\]
We obtain that in the case of rational $\alpha=p/q$ the
asymptotics of the eigenvalue distribution function $N_h(\lambda)$
as $h\rightarrow 0$ is given by
\[
 N_{h}(\lambda ) =h^{-1}
\sum_{k:\frac{4\pi^2}{p^2+q^2}k^2<\lambda}
\frac{1}{\pi\sqrt{p^2+q^2} }(\lambda -
\frac{4\pi^2}{p^2+q^2}k^2)^{1/2}+o(h^{-1}).
\]

Suppose that $\alpha\not\in\mathbb{Q}$. Let us the heat operator
$e^{-t\Delta_h}$  defined by the Laplace operator $\Delta_h$ in
${\mathbb R}^2$. For any $u\in L^2({\mathbb R}^2)$ the function
$u(t)=e^{-t\Delta_h}u$ is a solution of the heat equation
\[
\frac{\partial u}{\partial t}= -\Delta_hu, \quad t>0,
\]
with initial conditions $u(0)=u$. It can be shown that the kernel of
$e^{-t\Delta_h}$ on ${\mathbb R}^2$ has the form
\begin{multline*}
H_t(x,y,x_1,y_1)= \\ =\frac{h^{-1}}{4\pi
t}\exp\left[{-\frac{(x-x_1+\alpha(y-y_1))^2+h^{-2}(-\alpha
(x-x_1)+ y-y_1)^2}{4t(1+\alpha^2)}}\right].
\end{multline*}
Therefore, the kernel of $e^{-t\Delta_h}$  in ${\mathbb
T}^2$ is given by
\begin{multline*}
H_t(x,y,x_1,y_1)=\frac{h^{-1}}{4\pi t}
\sum_{k,l=-\infty}^{+\infty} \exp
[-\frac{(x-x_1+k+\alpha(y-y_1+l))^2}{4t(1+\alpha^2)}-\\
-\frac{(-\alpha (x-x_1+k)+ y-y_1+l)^2}{4th^{2}(1+\alpha^2)}].
\end{multline*}
Accordingly, the trace of $tr(e^{-t\Delta_h})$ in
$L^2({\mathbb T}^2)$ is given by
\[
tr(e^{-t\Delta_h})=\frac{h^{-1}}{4\pi
t}\sum_{k,l=-\infty}^{+\infty}{e^{-\frac{(k+\alpha l)^2}{4
t(1+\alpha^2)}-\frac{(-\alpha k+l)^2}{4th^2(1+\alpha^2)}}}.
\]
Since the series converges uniformly on $h\in (0,1]$ and, for any
$(k,l)\not=(0,0)$
\[
e^{-\frac{(k+\alpha l)^2}{4 t(1+\alpha^2)}-\frac{(-\alpha
k+l)^2}{4th^2(1+\alpha^2)}}\longrightarrow 0, \quad h\rightarrow
+0,
\]
we get
\[
{\rm tr}\, e^{-t\Delta_h}\longrightarrow \frac{h^{-1}}{4\pi
t},\quad h\rightarrow +0.
\]

From the other side, we have $$ {\rm tr}\, e^{-t\Delta_{h}}
=\int^{+\infty}_{-\infty} e^{-\lambda t}\,
d_{\lambda}N_{h}(\lambda ), $$ therefore, the proof  of
Theorem~\ref{th:main}  is  completed by the next lemma
(\cite[Lemma 5.2]{Sh:density}, cf. also \cite{av-simon}).

We say that a function $F$ on the real line  is  a distribution
function if $F$ is a non-decreasing left-continuous function such
that
\[
\lim_{\lambda \rightarrow -\infty} F(\lambda ) = 0.
\]
For any distribution function $F(\lambda )$, we denote by
$\tilde{F}(t)$ the Laplace transform of this function.

\begin{lem}
\label{Lemma 5.1} Let $F_{n}$ be a sequence of distribution
functions  on the real line such that:

(a) $F_{n}(\lambda ) = 0,\lambda  \leq a$, with the constant $a$,
not depending on $n$;

(b) $\vert F_{n}(\lambda )\vert  \leq  C e^{-\varepsilon \lambda
}, \lambda\in {\mathbb R}$, with the constants $C$ and
$\varepsilon $, not depending on $n$;

(c) there exists $\lim_{n \rightarrow \infty} \tilde{F}_{n}(t) =
\tilde{F}(t)$ for any $t>0$, where $F(t)$ is  a distribution
function.

Then $\lim_{n \rightarrow \infty} F_{n}(\lambda ) = F(\lambda )$
at  all  points  of  continuity  of the function $F(\lambda )$.
\end{lem}

\end{document}